\documentclass[12pt]{amsart}
\usepackage{amscd, amssymb}

\newtheorem{pr}{Proposition}
\newtheorem{lm}{Lemma}
\newtheorem{tm}{Theorem}

\newtheorem{pb}{Problem}

\newcommand\atsign{@}
\def\fra#1#2{\genfrac{}{}{0pt}{}{#1}{#2}}

\newcommand{\smallcup}{\mathbin{\text{\scriptsize$\cup$}}}
\newcommand\fsr{{\sc 5sr}}
\newcommand\tsr{{\sc 3sr}}
\newcommand\codim{\mathop{\rm codim}\nolimits}
\newcommand\Hom{\mathop{\rm Hom}\nolimits}

\newcommand\rk{\mathop{\rm rk}\nolimits}
\newcommand\Sym{\mathop{\rm Sym}\nolimits}
\newcommand{\Q}{{\mathbb Q}}
\newcommand{\Z}{{\mathbb Z}}
\newcommand{\R}{{\mathbb R}}
\newcommand{\C}{{\mathbb C}}

\newcommand{\bP}{{\mathbb P}}
\def\ga{\Gamma}
\def\feyn#1#2#3#4{\displaystyle\Bigl(
\fra{\displaystyle\hfill #1}{\displaystyle\hfill #2}{\rangle\hskip-3pt}
\frac{\ \ }{\ \ }{\hskip-3pt\langle}
\fra{\displaystyle #4\hfill}{\displaystyle #3\hfill}\Bigr)}
\def\dfeyn#1#2#3#4#5{\feyn{#1}{#2}{#3}{#4}^{#5}}

\def\fivearu#1#2#3#4#5{\feyn{#1\,#2{}^{^\nearrow}\hskip-8pt}{#3}{#4}{#5}}
\def\fiveard#1#2#3#4#5{\feyn{#1\hskip-3pt_{_\swarrow} #2}{#3}{#4}{#5}}
\def\fiveblu#1#2#3#4#5{\feyn{#1}{^{^\nwarrow}#2\,#3}{#4}{#5}}

\def\fivebrd#1#2#3#4#5{\feyn{#1}{#2\,#3{}_{_\searrow}\hskip-8pt}{#4}{#5}}
\def\fiveclu#1#2#3#4#5{\feyn{#1}{#2}{#3{}^{^\nearrow}\hskip-5pt #4}{#5}}

\def\fivecru#1#2#3#4#5{\feyn{#1}{#2}{#3\,#4{}^{^\nearrow}}{#5}}

\def\fivedlu#1#2#3#4#5{\feyn{#1}{#2}{#3}{\hskip-8pt^{^\nwarrow}#4\,#5}}

\def\dfiveard#1#2#3#4#5#6{\fiveard{#1}{#2}{#3}{#4}{#5}^{#6}}

\def\three#1#2#3#4#5{\feyn{#1}{}{#2}{#3,#4;#5}}
\def\dthree#1#2#3#4#5#6{\three{#1}{#2}{#3}{#4}{#5}^{#6}}

\begin{document}
\title{Associativity relations in quantum cohomology}
\author{ Andrew Kresch$^1$}
\date{11 March 1997}
\begin{abstract}
We describe interdependencies among the
quantum cohomology associativity relations.
We strengthen the first reconstruction theorem of
Kontsevich and Manin by identifying a subcollection
of the associativity relations which
implies the full system of WDVV equations.
This provides a tool for identifying non-geometric
solutions to WDVV.
\end{abstract}
\maketitle
\pagestyle{plain}
\footnotetext[1]{Funded by a
Fannie and John Hertz Foundation Fellowship
for graduate study.}
\section{Introduction}
\label{intro}

The geometry of moduli spaces of stable maps of genus 0
curves into a complex projective manifold $X$ leads to
a system of quadratic equations in the tree-level
(genus 0) {\em Gromov-Witten numbers} of $X$.

These quadratic equations were written down by physicists before the
geometry was worked out rigorously.
In all the worked examples, the equations
were seen to determine a solution, uniquely and consistently,
from starting data.
The beautiful paper of Di Francesco and Itzykson, \cite{a},
presents a number of examples in this context.

In one of the foundational papers in the area of quantum cohomology,
\cite{b}, the authors remark on the overdeterminedness of the
system of equations, and pose the
question whether the seemingly redundant equations
follow algebraically from the useful ones.

This same paper presents the
{\em first reconstruction theorem}, which applies to
manifolds $X$ such that $H^*(X,\Q)$ is generated by $H^2(X)$.
This result gives an effective procedure for solving for
genus 0 Gromov-Witten numbers from starting data using the
quadratic relations
(since this entire paper is concerned only with the genus 0 invariants,
we omit explicit mention of genus from now on).

This paper continues in the spirit of these early investigations
into the structure of these relations.
After a review in which notation is introduced and the
basic problem is set up
(section \ref{review}), some results
are presented (sections \ref{threesi} and \ref{fivesi})
explicitly showing interdependencies among these relations.

Then follows a generalization of the first reconstruction
theorem.
In this, we are forced to keep the hypothesis that the
cohomology ring of $X$ be generated by divisors.
We then show that an initial collection of numbers and
relations determines, purely algebraically, the entire
system of relations ({\em strong reconstruction})
as well as the Gromov-Witten numbers.

Finally, some examples are worked out.
The structure of the associativity relations in the case $X$ has
dimension 2 is particularly easy to describe.
Section \ref{srttwo} gives an
explicit picture of strong reconstruction in this case.
Section \ref{examples} presents
some higher-dimensional examples.

Focusing on the algebra, rather than the geometry, of
quantum cohomology leads to some generalizations.
First, we can do away with $X$ entirely and work just with
its cohomology ring, or more generally any
{\em positively graded Gorenstein $\Q$-algebra}
with socle in degree $n\ge2$
($n$ plays the role of the dimension of $X$,
and we also must assume that the zeroth graded piece is isomorphic to $\Q$).
The canonical class $K_X$ plays an important role in
geometry since it determines the expected dimension of the
moduli spaces, but is less important to us:
we can change $K_X$ at will, or drop it entirely from
the discussion.
Finally, the very equations lead us to solutions which do
not come from geometry.
Some of the examples exhibit that the geometric solution,
with its family of obvious rescalings, may live in an even larger
solution space.

In an appendix (section \ref{geoground}),
we sketch some of the geometry that the rest of the paper ignores.
This is used to motivate the result in section \ref{fivesi},
although this result is purely algebraic.
This is in no sense a substitute for a survey
of quantum cohomology; we recommend \cite{c} as
starting point to reader unacquainted with the subject.

The author, who spent the 1996--97 year at the Mittag-Leffler Institute,
would like to take this opportunity to thank the staff and
organizers there for providing a comfortable and stimulating
atmosphere for research.
The author would also like to
express his appreciation to P. Belorousski, C. Faber, W. Fulton, T. Graber,
B. Kim, S. Kleiman, and R. Pandharipande for insightful conversation.

\section{The basic problem}
\label{review}
Setting up and solving the system of
quadratic equations in the Gromov-Witten numbers requires
a ring which, in an appropriate sense, looks like a
cohomology ring of a manifold.

Precisely, let $A$ be a
positively graded Gorenstein $\Q$-algebra with socle
in degree $n\ge2$
such that the zeroth graded piece of $A$ is isomorphic to $\Q$,
and assume a nonzero choice of
$\varphi\in\Hom_A(\Q,A)\simeq\Q$ is made.
We call the degree of an element its {\em codimension},
and denote the $k^{\rm th}$ graded piece by $A^k$.
Then $\varphi$ is an isomorphism of $\Q$ with $A^n$.
Denote by $\int$ the composite of projection with this isomorphism:
$$A\to A^n\stackrel{\varphi^{-1}}\to\Q.$$
If $\{T_i\,|\,i\in I\}$ is a basis for $A$ as a $\Q$-module,
where $I$ is a (necessarily finite) indexing set,
let $g_{ij}=\int T_i\cdot T_j$.
This is a nondegenerate pairing; we denote by $(g^{ij})$ the inverse matrix.
We make the further assumption that $A^1$ is nonzero.

In order to get a well-defined system of equations, we
need two more pieces of data: a maximal-rank integral lattice $\Lambda$
and a strongly convex polyhedral cone $\Theta$,
both contained in ${(A^1)}^*\otimes\R$.

Let $B=\{T_1,T_{\sigma_1},\ldots,T_{\sigma_r},
T_{\tau_1},\ldots,T_{\tau_s}\}$ be
a basis of $A$ as a $\Q$-module, consisting of homogeneous elements, with
$T_1=1$,
$T_{\sigma_i}\in A^1$ for each $i$,
$\codim T_{\tau_j}\ge2$ for each $j$,
and $\codim T_{\tau_j}\le\codim T_{\tau_k}$ whenever $j\le k$.
This defines $r$ and $s$ as the ranks of the first
and higher graded pieces of $A$, respectively; recall also that $n$ is
the top grading, necessarily equal to $\codim T_{\tau_s}$.
Assume the choice of basis is made so that $\Lambda$ is
dual to the integral span of $T_{\sigma_1},\ldots,T_{\sigma_r}$.

Geometry provides many examples of such data:
$A^*=H^{2{*}}(X,\Q)$, where $X$
is a complex projective manifold with homology
only in even dimensions, $\int=\int_X$,
$\Lambda=H_2(X,\Z)/{\rm torsion}$,
$\Theta=$ dual to ample cone.

Let $C=\Theta\cap\Lambda\setminus\{0\}$.
We define the set of unknown numbers to be the collection of all
$N(\beta;d_1,\ldots,d_s)$
with $\beta\in C$ and $d_j\ge0$ for all $j$.
For any $\beta\in C$,
denote by $c_i$ the pairing of $\beta$ with $T_{\sigma_i}$.
We now define a system of equations in these unknowns, one for
each 4-tuple $(i,j,k,l)$ of elements of the basis indexing set
$I=\{1,\sigma_1,\ldots,\sigma_r,\tau_1,\ldots,\tau_s\}$
and each {\em degree} $(\beta;d_1,\ldots,d_s)$.

To do this, introduce the {\em potential function}, a formal
Laurent series
$$\Phi=\Phi^{\rm cl} + \Gamma$$
in variables $y_i$, $i\in I$,
composed of the {\em classical part}
$$\Phi^{\rm cl}=\frac{1}{6}\sum_{i,j,k\in I}
\left(\int T_i\cdot T_j\cdot T_k\right) y_iy_jy_k$$
and the {\em quantum correction}
$$\Gamma=\sum_{\fra{\scriptstyle\beta\in C}{\scriptstyle d_1,\ldots,d_s\ge0}}
N(\beta;d_1,\ldots,d_s)
e^{c_1y_{\sigma_1}}\cdots e^{c_ry_{\sigma_r}}
\frac{y_{\tau_1}^{d_1}}{d_1!}\cdots\frac{y_{\tau_s}^{d_s}}{d_s!}.$$
Then any $(i,j,k,l)$ determines a differential equation
\begin{equation}
\label{assrel}
\sum_{e,f}\Phi_{ije}g^{ef}\Phi_{fkl}=
\sum_{e,f}\Phi_{jke}g^{ef}\Phi_{fil},
\end{equation}
where we have used subscripts to denote partial differentiation.
Isolating the coefficient of
$e^{c_1y_{\sigma_1}}\cdots e^{c_ry_{\sigma_r}}
y_{\tau_1}^{d_1}\cdots y_{\tau_s}^{d_s}$
on each side produces a quadratic equation in $N$'s,
which we call an {\em associativity relation}
(they imply associativity of the so-called quantum product; see \cite{c}).
We follow Dubrovin in calling the system of
equations (\ref{assrel}) the {\em WDVV equations}
(after E. Witten, R. Dijkgraaf, H. Verlinde and E. Verlinde).
A particular WDVV equation is represented symbolically by
an equivalence of Feynman diagrams
$$\feyn ijkl \ \sim\ \Bigl(
\begin{picture}(40,10)(-20,17)
\put(0,16){\line(0,1){8}}
\put(0,16){\line(3,-1){10}}
\put(0,16){\line(-3,-1){10}}
\put(0,24){\line(3,1){10}}
\put(0,24){\line(-3,1){10}}
\put(12,8){\em k}
\put(-17,8){\em j}
\put(12,26){\em l}
\put(-17,26){\em i}
\end{picture}\Bigr)$$
We adopt the notation $\feyn ijkl$ to refer to the WDVV equation
(\ref{assrel}) and
$\dfeyn ijkl{(\beta;d_1,\ldots,d_s)}$ to refer to
a particular associativity relation.

We would like to rewrite (\ref{assrel}) in which we split $\Phi$
into its classical and quantum parts.
If $T_i\cdot T_j=\sum_h q_h T_h$ then denote 
$\sum_h q_h\ga_{hkl}$ by $\ga_{(ij)kl}$.
Then by rewriting (\ref{assrel}) we get
\begin{equation}
\label{altassrel}
\ga_{ij(kl)}+\ga_{(ij)kl}-\ga_{jk(il)}-\ga_{(jk)il}
= \sum_{e,f}\ga_{jke}g^{ef}\ga_{fil}-
\sum_{e,f}\ga_{ije}g^{ef}\ga_{fkl}.
\end{equation}
Thus we have split the WDVV equation into the
{\em linear contribution} (left-hand side) and
{\em quadratic contribution} (right-hand side).

We emphasize that $\dfeyn ijkl{(\beta;d_1,\ldots,d_s)}$ is
an equation; it can be scaled or combined linearly with
other equations, and it can imply other equations.
We will sometimes need to distinguish the linear and quadratic contributions
to a particular associativity relation, each being the polynomial
appearing on the appropriate side of (\ref{altassrel}).

With this set-up, we can now state
\begin{pb}
\label{mainprob}
Given $A$, $\int$, $\Lambda$, $\Theta$ as above,
find solutions in rational numbers $N(\beta;d_1,\ldots,d_s)$
to the full set of WDVV equations (\ref{assrel}).
\end{pb}

For the purposes of ordering elements of $C$, we let
$\omega$ be an element of the interior of the dual cone to
$\Theta$ (which must be nonempty since $\Theta$ is
strongly convex).
Then
$\langle\beta,\omega\rangle > 0$ for all $\beta\in C$,
and for each $k$, the set
$\{\beta\in\Z\langle\sigma^*_1,\ldots,\sigma^*_r\rangle\,|\,
\langle\beta,\omega\rangle < k\}$ is finite.

We must remark on four distinct uses of the Feynman symbols.
First, with $i,j,k,l\in I$, $\feyn ijkl$ is as above.
Second, for $\xi, \pi, \rho, \sigma\in A$,
$\feyn\xi\pi\rho\sigma$ refers to the equation
obtained by writing each element in
terms of the basis and summing in a multilinear fashion.
Third, for subsets $\Xi, \Pi, {\rm P}, \Sigma$ of $A$,
$\feyn \Xi\Pi{\rm P}\Sigma$ refers to the collection of
all $\feyn\xi\pi\rho\sigma$ with $\xi$ in $\Xi$, etc.
Finally, as a special case, for integers $w, x, y, z$,
$\feyn wxyz$ refers to $\feyn {A^w}{A^x}{A^y}{A^z}$.
Also, if $i,j,k,l,m\in I$, we write $\feyn{ij}klm$
as shorthand for
$\feyn{T_i\cdot T_j}{T_k}{T_l}{T_m}$.

\section{How many relations?}
\label{howmany}

As we have seen, to each choice of $i,j,k,l\in I$ there corresponds a
WDVV equation.
Thus, if the vector space $A$ has dimension $r$ then
the number of distinct WDVV equations is of the order of
magnitude $r^4$.

It is an exercise in combinatorics to provide a precise count.
First, if any of $i$, $j$, $k$, $l$ is 1
(i.e., indexes the identity element of $A$),
then the corresponding WDVV equation is a trivial identity.
If $k=i$ or $l=j$ then the equation is trivial.

There are symmetries.
If we swap $i$ and $j$ and swap $k$ and $l$ then the
WDVV equation remains unchanged.
Swapping two symbols on a diagonal of the Feynman diagram
only changes the WDVV equation by a sign.
If $i$, $j$, $k$, and $l$
are all distinct, then the three relations one obtains by
cyclically permuting $i$, $j$, $k$ are linearly related:
\begin{equation}
\label{twooutofthree}
\feyn ijkl + \feyn jkil + \feyn kijl = 0.
\end{equation}

A tally shows that the number of distinct nontrivial WDVV equations
(modulo sign) is
$$\frac{r^4-6r^3+15r^2-18r+8}{8},$$
while if we count only two out of the three distinct WDVV equations
involving four distinct symbols --- since by (\ref{twooutofthree})
any two imply the third,
a fact we refer to as the {\em two-out-of-three} implication --- the
count is
$$\frac{r^4-4r^3+5r^2-2r}{12}.$$
Since there is only (at most) one unknown number in each degree,
the count above
is, na\"{\i}vely, the factor of overdeterminedness in
the system of associativity relations.

\section{The three symbols relation}
\label{threesi}

Let $i,j,k,l,m\in I$.
Let $\Phi=\Phi^{\rm cl} + \ga$ be the potential function.
The following algebraic identity,
called the {\em three symbols identity}, holds:
\begin{eqnarray}
\lefteqn{\frac{\partial}{\partial y_m}
\Bigl(\sum_{e,f} \Phi_{ije}g^{ef}\Phi_{fkl} -
 \sum_{e,f} \Phi_{jke}g^{ef}\Phi_{fil}\Bigr) + {}} \nonumber \\
& & \frac{\partial}{\partial y_j}
\Bigl(\sum_{e,f} \Phi_{ile}g^{ef}\Phi_{fkm} -
 \sum_{e,f} \Phi_{lke}g^{ef}\Phi_{fim}\Bigr) + {} \label{threesymbol}\\
& & \frac{\partial}{\partial y_l}
\Bigl(\sum_{e,f} \Phi_{ime}g^{ef}\Phi_{fkj} -
 \sum_{e,f} \Phi_{mke}g^{ef}\Phi_{fij}\Bigr) = 0. \nonumber
\end{eqnarray}

Let $(\beta,d)$ be a degree.
Define $e_{\sigma_i}=0$ and $e_{\tau_j}=(0,\ldots,1,\ldots,0)$
with the 1 in the $j^{\rm th}$ place.
Then (\ref{threesymbol}) gives us
\begin{pr}
\label{threesymrel}
Suppose $i,j,k,l,m\in I$ with $\codim T_m\ge2$,
and let $(\beta,d)$ be a degree with
$d_m\ge1$.
Then relations $\dfeyn ilkm{(\beta;d+e_j-e_m)}$ and
$\dfeyn imkj{(\beta;d+e_l-e_m)}$ together
imply $\dfeyn ijkl{(\beta;d)}$.
\end{pr}
This we call the {\em three symbols relation} (\tsr), and denote by
the diagram $\dthree ikjlm{(\beta;d)}$.

We now record one application of \tsr.

\begin{lm}
\label{lemma}
With the notation of Problem \ref{mainprob},
suppose $(\beta;d)$ is a degree with $d\ne0$.
Then the collection of all relations in degrees $(\beta;d')$ with
$\sum_id'_i=\left(\sum_id_i\right)-1$ implies
$\dfeyn A{A^1}A{A^1} {(\beta;d)}$.
\end{lm}

Indeed, if $d_m\ge1$ then
$\three ij11m$ yields $\dfeyn i1j1{(\beta;d)}$
for any $i,j\in I$.

\section{The five symbols relation}
\label{fivesi}

Given $i,j,k,l,m\in I$, the following algebraic identity holds:
\begin{equation}
\label{mdiag}
\sum_{e,f} \ga_{ij(me)}g^{ef}\ga_{fkl} =
\sum_{e,f} \ga_{kl(me)}g^{ef}\ga_{fij}.
\end{equation}
Recall, if $T_m\cdot T_e=\sum_q t_q T_q$ then
by $\ga_{ij(me)}$ we mean $\sum_q t_q\ga_{ijq}$.
Now (\ref{mdiag}) follows by observing that with
$g_{abc}=\int T_a\cdot T_b\cdot T_c$
we have $t_q=\sum_p g_{mep} g^{pq}$, and now the coefficient
$\sum_{e,p} g^{ef}g_{mep}g^{pq}$ of
$\ga_{ijq}\ga_{fkl}$ on the left-hand side is symmetric in $f$ and $q$.

The appendix gives geometric motivation for (\ref{mdiag}).

We write the expression 
$\sum \ga_{ij(me)}g^{ef}\ga_{fkl} -
\sum \ga_{kl(me)}g^{ef}\ga_{fij}$ and
add to it the four additional expressions obtained by permuting
the variables $i,j,k,l,m$ cyclically.
We use the identity coming from the associativity relation
$\feyn emij$
and its cyclic translates to obtain
\begin{eqnarray*}
\lefteqn{0 =
\ga_{ij(me)}\ga_{fkl} + \ga_{jk(ie)}\ga_{flm} + \ga_{kl(je)}\ga_{fmi} +
\ga_{lm(ke)}\ga_{fij} + \ga_{mi(le)}\ga_{fjk} } \\
 & & {} - \ga_{mi(je)}\ga_{fkl} - \ga_{ij(ke)}\ga_{flm} -
 \ga_{jk(le)}\ga_{fmi} -
\ga_{kl(me)}\ga_{fij} - \ga_{lm(ie)}\ga_{fjk} \\
\lefteqn{\phantom0 =
\ga_{(mi)je}\ga_{fkl} + \ga_{(ij)ke}\ga_{flm} + \ga_{(jk)le}\ga_{fmi} +
\ga_{(kl)me}\ga_{fij} + \ga_{(lm)ie}\ga_{fjk} } \\
 & & {} - \ga_{(ij)me}\ga_{fkl} - \ga_{(jk)ie}\ga_{flm} -
 \ga_{(kl)je}\ga_{fmi} -
\ga_{(lm)ke}\ga_{fij} - \ga_{(mi)le}\ga_{fjk}.
\end{eqnarray*}
We have omitted summations symbols and $g^{ef}$'s to save space.
We have also omitted the (cubic) terms obtained by substituting the
quadratic contributions of the associativity relations, but
the key point is that these cancel.

The final expression above is the quadratic contribution of a sum of
associativity relations,
conveniently written
\begin{equation}
\label{fivearound}
\feyn{mi}jkl - \feyn m{ij}kl + \feyn mi{jk}l -
\feyn mij{kl} + \feyn{lm}ijk.
\end{equation}

Since the linear contribution of (\ref{fivearound}) vanishes, as may
be checked, we have, at least formally, that the indicated associativity
relations imply the vanishing of (\ref{fivearound}).

We turn this into a precise, practical statement by grading the
terms in $\ga$ by degree.
We use the notation of Problem \ref{mainprob}.
Rewrite the above, isolate the coefficient of
some degree $(\beta;d)$, and note that then every quadratic
term is a sum over
$\beta_1+\beta_2=\beta$ with
$\langle \beta_i,\omega \rangle < \langle \beta,\omega \rangle$ for $i=1,2$.
This establishes

\begin{pr}
\label{fivesymrel}
Suppose $i,j,k,l,m\in I$, and let $(\beta;d)$ be a degree.
The collection of relations consisting of
$\dfeyn ijke{(\beta';d')}$ and its cyclic translates through $\{i,j,k,l,m\},$
for all $e\in I$ and all degrees $(\beta',d')$ with
$\langle \beta',\omega\rangle < \langle \beta,\omega\rangle$
and $d'\le d$ (componentwise),
implies the relation
\begin{eqnarray*}
\lefteqn{\dfeyn{mi}jkl{(\beta;d)} - \dfeyn m{ij}kl{(\beta;d)} +
\dfeyn mi{jk}l{(\beta;d)}}\hspace{100pt} \\
& & {} - \dfeyn mij{kl}{(\beta;d)} + \dfeyn{lm}ijk{(\beta;d)}.
\end{eqnarray*}
\end{pr}

We call this the {\em five symbols relation} (\fsr).
We employ the notation $\dfiveard mijkl{(\beta;d)}$ to describe the above
relation.

\section{Strong reconstruction for $n=2$}
\label{srttwo}
To illustrate an application of the three and five symbols relations,
we work out a strong reconstruction theorem for $n=2$,
where the associativity relations are simple to organize.
When $n=2$, there is only one $\tau$, and there are three
types of associativity relations:
\begin{enumerate}
\item[(i)] $\feyn \tau \tau {\sigma_j} {\sigma_i}$;
\item[(ii)] $\feyn \tau {\sigma_i} {\sigma_j} {\sigma_k}$;
\item[(iii)] $\feyn {\sigma_i} {\sigma_j} {\sigma_k} {\sigma_l}$.
\end{enumerate}

Assume for simplicity that, as in the geometric situation,
there is a canonical class $K$ which dictates that there
is at most one nontrivial number $N(\beta;d)$ in each
curve class $\beta$, namely when
$d=\langle \beta,-K \rangle - 1 \ge 0$
(for the general case, see the next section).
The potential function is composed of
\begin{eqnarray*}
\Phi^{\rm cl} &=& \frac{1}{2} y_1^2 y_\tau +
\frac{1}{2} \sum_{e,f=1}^r g_{ef} y_1 y_{\sigma_e} y_{\sigma_f}, \\
\ga &=& \sum_{\fra{\scriptstyle\langle\beta,-K\rangle\ge1}
{\fra{\scriptstyle\beta=\sum c_i\sigma^*_i\in C}
{\scriptstyle d=\langle \beta,-K \rangle - 1}}}
N(\beta;d)e^{c_1 y_{\sigma_1}}\cdots e^{c_r y_{\sigma_r}}
\frac{y_\tau^d}{d!}
\end{eqnarray*}
(we write $g_{ef}$ for $g_{\sigma_e\sigma_f}$).

Suppose we are given all $N(\beta;d)$ with $d\le2$, and suppose
that these satisfy
$$\dfeyn \tau {\sigma_i} {\sigma_j} {\sigma_k} {(\beta;0)}{\rm\ and\ }\dfeyn
{\sigma_i} {\sigma_j} {\sigma_k} {\sigma_l} {(\beta;0)}$$
for all $i,j,k,l\in\{1,\ldots,r\}$ and all $\beta$.
We claim that relations of type (i) allow us to solve for all further
$N(\beta;d)$ ({\em reconstruction}) and that the numbers thus
obtained satisfy the full system of WDVV equations
({\em strong reconstruction}).

Indeed, by the three symbols relation,
$$
\dfeyn \tau {\sigma_j} {\sigma_i} {\sigma_l} {(\beta;0)} {\rm\,and\,\ }
\dfeyn \tau {\sigma_j} {\sigma_k} {\sigma_l} {(\beta;0)}
\,\Rightarrow\,\dfeyn {\sigma_i} {\sigma_j} {\sigma_k} {\sigma_l} {(\beta;1)},
$$
and thus the hypothesis implies relations (ii) $(\beta;1)$.
Similarly, relations (i) $(\beta;d)$ imply
(ii) $(\beta;d+1)$ and (iii) $(\beta;d+2)$.

Inductively on $d$,
assume all $N(\beta;d')$ known and all relations satisfied for
$\langle \beta, -K \rangle < d+4$.
Relation $\dfeyn \tau \tau {\sigma_j} {\sigma_i} {(\beta;d)}$
reads $g_{ij}\ga_{zzz}^{(\beta;d)}=Q_{ij}^{(\beta;d)}$, where
$Q_{ij}^{(\beta;d)}$ is a quadratic expression in
known quantities.
Now $\fivebrd\tau {\sigma_k} {\sigma_l} {\sigma_j} {\sigma_i}$
tells us $g_{kl}Q_{ij}^{(\beta;d)} = g_{ij}Q_{kl}^{(\beta;d)}$,
which says we can solve for
$\ga_{zzz}^{(\beta;d)}$
(that is, $N(\beta;d+3)$) satisfying (i), and
we have just seen that the relations of type (i) imply all the relations
in degree $\beta$.

\section{Strong reconstruction theorem}
\label{srt}
The main result here is a generalization of the
first reconstruction theorem of Kontsevich and Manin \cite{b}.
Working in the same class of cohomology rings, namely those
generated by elements in codimension 1, we prove that an identifyable
collection of numbers satisfying an identifyable collection of
relations gives us strong reconstruction, i.e.\ allows us
to solve uniquely for numbers
satisfying the complete system of associativity relations.
In case $-K$ is ample, we need only a finite collection
of numbers and relations as starting data.

The result is
\begin{tm}
\label{strongr}
With the notation of Problem \ref{mainprob},
suppose $A$ is generated by $A^1$.
Then the collection of $N(\beta;d)$ with $\sum_{i=1}^s d_i\le2$
extends to a solution to WDVV if and only if
$\dfeyn A{A^1}A{A^1} {(\beta;0)}$
is satisfied for all $\beta$.
\end{tm}

We begin by organizing notation.
By hypothesis,
we may assume the basis $B$ chosen such that
for each $j$, $1\le j\le s$, there exists $i_j\in\{1,\ldots,r\}$
and $\mu_j\in I$ such that $T_{\tau_j}=T_{\sigma_{i_j}}\cdot T_{\mu_j}$.

We wish to impose a partial order on the collection of degrees
$d=(d_1,\ldots,d_s)$ with fixed $|d|:=\sum_{i=1}^s d_s$, such that
$(d_1,\ldots,d_s)$ precedes $(d_1,\ldots,d_i+1,\ldots,d_j-1,\ldots,d_s)$
for any $i < j$.
A convenient way is to order by $\sum id_i$.

Let us give an outline of the proof of the theorem.
Inducting on $\langle \beta, \omega \rangle$, then on $|d|$,
then downwards on $\sum id_i$,
we verify all associativity relations in degree $(\beta;d)$,
showing that those of the form
$$\dfeyn {\mu_j}{\sigma_{i_j}}{\tau_k}{\tau_l}{(\beta;d)}$$
with
$\codim T_{\tau_j}\le\codim T_{\tau_k}\le\codim T_{\tau_l}$
and $\max(j,k,l)\le\min\{m\,|\,d_m\ne0\}$
determine the numbers
$N(\beta;d+e_j+e_k+e_l)$ (here $e_i=(0,\ldots,1,\ldots,0)$ with
1 in the $i^{\rm th}$ position).

\section{Proof of the strong reconstruction theorem}
\label{srtproof}
The induction breaks up into an outer induction on degrees
and an inner induction within each degree.
The outer induction proceeds
with respect to the partial order: $(\beta',d')\prec(\beta,d)$
if $\langle \beta',\omega\rangle < \langle \beta,\omega\rangle$
and $|d'|\le|d|$;
$\beta'=\beta$ and $|d'|<|d|$;
or $\beta'=\beta$, $|d'|=|d|$, and
$\sum i d'_i > \sum i d_i$.
The inner induction is on $(u,c,a,b)$ with $u$
(corresponding to $\codim T_{\mu_j}$ above) up from 1,
$c$ ($=\codim T_{\tau_k}+\codim T_{\tau_l}$) up from $2(u+1)$, and
$a$ ($=\codim T_{\tau_k}$) up from $u+1$ to $[c/2]$.
Define $b=c-a$; then we always have $a\le b$.

The induction hypothesis, at a given step $(\beta,d,u,c,a,b)$,
consists of all relations in previous degrees plus all
numbers they refer to (i.e., all $N(\beta';d'+e')$
with $(\beta',d')\prec(\beta,d)$, $|e'|\le3$),
plus, in the current degree,
all $\feyn z1xy$ with $\min(x,y,z)<u$,
all $\feyn u1xy$ and $\feyn x1yu$
with $x\ge u+1$, $y\ge u+1$,
and either $x+y<c$ or $x+y=c$ with $\min(x,y)<a$,
together with all numbers these relations refer to.

In any degree, for any
integers $x$ and $y$, $\feyn x1y1$ follows either by hypothesis ($d=0$)
or by the induction hypothesis and Lemma \ref{lemma} ($|d|\ge1$).
When $u\ge2$ we obtain $\feyn x1yu$ for $x\ge u$ and
$y\ge u$ from $\fivedlu x1y{u{-}1}1$, and now
$\feyn u1ux$ for $x>u$ from
$\fivearu {u{-}1}11ux$.

The main step is to deduce $\feyn u1ab$.
Here the linear terms coming from the associativity relation
possibly involve new $N$'s.
We divide this into two steps.

First, we show it suffices to prove a distinguished set
of $\feyn u1ab$.
Let $S$ be the set of relations
$\feyn {\mu_j}{\sigma_{i_j}}{\tau_k}{\tau_l}$ with
$\codim T_{\tau_j}=u+1$,
$\codim T_{\tau_k}=a$, and $\codim T_{\tau_l}=b$.
We claim that $S$ (and the new $N$'s referred to)
implies $\feyn u1ab$.
Indeed, if
$\codim T_\mu=u$ and $\codim T_\sigma=1$ with
$T_\sigma\cdot T_\mu=\sum \lambda_j T_{\tau_j}$, then
comparing $\fiveclu \mu\sigma{\sigma_{i_k}}{\mu_k}{\tau_l}$ with
$\sum \lambda_j \fiveclu {\mu_j}{\sigma_{i_j}}{\sigma_{i_k}}{\mu_k}{\tau_l}$
establishes $\feyn \mu\sigma{\tau_k}{\tau_l}$ from the
relations in $S$.

For the second step, we establish all relations in $S$.
Each $\feyn {\mu_j}{\sigma_{i_j}}{\tau_k}{\tau_l}$ in $S$
involves the variable $N(\beta;d+e_j+e_k+e_l)$.
For $a$, $b$, $u+1$ distinct, there is a one-to-one correspondence
between elements of $S$ and such variables.
In other cases, we shall need symmetrizing arguments to show any two
relations in $S$ sharing a common such variable are equivalent.
In case $a=b$, the two-out-of-three implication gives
$\feyn {\mu_j}{\sigma_{i_j}}{\tau_k}{\tau_l}\Leftrightarrow
\feyn {\mu_j}{\sigma_{i_j}}{\tau_l}{\tau_k}$.
In case $a=u+1$, we get
$\feyn {\mu_j}{\sigma_{i_j}}{\tau_k}{\tau_l}\Leftrightarrow
\feyn {\sigma_{i_j}}{\mu_j}{\tau_k}{\tau_l}\Leftrightarrow
\feyn {\mu_k}{\sigma_{i_k}}{\tau_j}{\tau_l}$
by two-out-of-three and
$\fivecru {\sigma_{i_j}}{\mu_j}{\sigma_{i_k}}{\mu_k}{\tau_l}$.

Thus, it suffices to establish only those
$\feyn{\mu_j}{\sigma_{i_j}}{\tau_k}{\tau_l}\in S$
such that $j\le k\le l$.
In case $d_m\ge1$ for some $m<l$,
$\three {\mu_j}{\tau_k}{\sigma_{i_j}}{\tau_l}{\tau_m}$
establishes $\feyn{\mu_j}{\sigma_{i_j}}{\tau_k}{\tau_l}$.
Otherwise, $N(\beta;d+e_j+e_k+e_l)$ is actually an unknown,
so solving $\feyn{\mu_j}{\sigma_{i_j}}{\tau_k}{\tau_l}$
establishes simultaneously the number and
the relation.
Finally, two-out-of-three
establishes $\feyn u1ba$ from
$\feyn u1ab$.

Having finished the inner induction,
to establish general $\feyn wxyz$
is an easy induction on $\min(w,x,y,z)$, using
\fsr{} by decomposing the entry of lowest codimension.

\section{Examples}
\label{examples}
Since we wish not to stray far from geometry, we focus mainly on rings
of the form $A=A^*_\Q X:=A^*X\otimes\Q$,
where $X$ is a projective manifold, and
with a dimension restriction on numbers coming from a class
$K\in A^1$:
$N(\beta;d)=0$ unless
$\sum_j d_j
(\codim T_{\tau_j}-1)=\langle\beta, -K\rangle + n - 3,$
where $n=\dim X$.
Unless otherwise stated, $K=K_X$, the canonical class on $X$.

{\bf Example 1.} $X=\bP^n$. The ring $A=A^*_\Q\bP^n$ has $\rk A^1=1$,
so strong reconstruction dictates a vacuous set of relations
(a Feynman diagram with a symbol appearing twice, in opposite corners,
determines a trivial relation),
and the solutions to WDVV correspond exactly to choices
of $N(\beta;d)$ with $|d|\le2$.
With $K=K_X$, there is only one such $N$, so there
are no solutions beyond the geometric solution and its rescalings.

Let $h=c_1({\mathcal O}(1))$.  When $K=-bh$ with $1\le b\le n-1$, there
is more than one $N(\beta;d)$ with $|d|\le2$, and hence
there is a family of solutions to WDVV of dimension $>1$.
Some of these solutions have geometric significance.
If $V\stackrel{i}\hookrightarrow\bP^{n+k}$ is a $(s_1,\ldots,s_k)$-complete
intersection with $b=n+k+1-\sum s_j\ge 1$, and $K=-bh$, then
$K_V=i^*K$, and we get a solution to WDVV for $A=A^*_\Q\bP^n$ with this
``wrong'' canonical class by setting
$N(\beta;d)$ equal to the sum over ${\beta'}\in H_2V$
with $i_*{\beta'}=\beta$ of the (geometric) Gromov-Witten number
for $V$
corresponding to $\bigotimes (i^*T_{\tau_j})^{\otimes d_j}$
in curve class $\beta'$ (cf.\ \cite{g}, sec.\ 4).

{\bf Example 2.} A toric manifold.
In the integral lattice of rank 3, let
$v_1=(1,0,0)$, $v_2=(0,1,0)$,
$v_3=(0,0,1)$, $v_4=(-2,-1,-1)$, $v_5=(-1,0,0)$,
and let the toric manifold $X$ be defined by the fan consisting
of the cones $\langle v_1,v_2,v_3\rangle $,
$\langle v_1,v_2,v_4\rangle$,
$\langle v_1,v_3,v_4\rangle$,
$\langle v_2,v_3,v_5\rangle$,
$\langle v_2,v_4,v_5\rangle$,
$\langle v_3,v_4,v_5\rangle$.
If $D_i$ denotes the divisor corresponding to vector $v_i$,
then $D_1$, $D_2$ span $A^1X$ with the dual to the ample cone
generated by the dual basis elements $D^*_1$, $D^*_2$.
We have $A=A_\Q^*X=\Q[D_1,D_2]/(D_1^2-2D_1D_2, D_2^3)$
with $\Q$-basis
$\{ 1, D_1, D_2, D_1 D_2, D_2^2, D_1D_2^2 \}$
and $-K=2 D_1 + D_2$.

Given this setup, the strong reconstruction theorem dictates
21 relations involving 17 variables.
When written out, this system of equations reduces to
\begin{eqnarray*}
N(D^*_1;2,0,0) &\hskip-4pt=\hskip-4pt& N(D^*_1;0,0,1); \\
N(D^*_1+D^*_2;0,1,1) &\hskip-4pt=\hskip-4pt& - N(D^*_1;0,0,1)
N(D^*_2;0,1,0); \\
0 &\hskip-4pt=\hskip-4pt&
N(D^*_1;0,0,1)\left[N(2D^*_2;0,2,0)+N(D^*_2;0,1,0)^2\right]
\end{eqnarray*}
(with the rest of the variables all zero).
Thus the solution space to WDVV has two irreducible components: one, containing
the geometric solution, with its expected two-dimensional family
of rescalings, and another, supported in curve classes lying
along one edge of the dual
to the ample cone, with a one-dimensional family of rescalings plus
dependence on another free parameter.

{\bf Example 3.} $X=G(2,4)$, $X=\Sym^2\bP^2$.  Both have isomorphic
cohomology rings (up to scale factors), not generated by divisors,
so we are outside the scope of the strong reconstruction
theorem.  We illustrate here a technique which allows one to prove,
subject to a genericity hypothesis on starting data, a strong
reconstruction result in this case.
The starting data according to the statement of strong
reconstruction consists of just one number
when $X=G(2,4)$ and 3 numbers
when $X=\Sym^2\bP^2$.
The result we seek is that modulo a genericity hypothesis
(some quantity in the starting data being nonzero) any
choice of starting data extends to a solution to WDVV.
It is helpful to recall the dimension condition
on relations.
For $\dfeyn \xi\pi\rho\sigma{(\beta;d)}$ to be nontrivial requires
\begin{eqnarray*}
\lefteqn{\langle\beta, -K\rangle -
\sum_{j=1}^s d_j(\codim T_{\tau_j}-1) = } \hspace{80pt} \\
& & \codim \xi + \codim \pi + \codim \rho + \codim \sigma - n.
\end{eqnarray*}

One may take as
cohomology basis the powers of the ample generator $h$ of $A^1X$,
plus an extra codimension 2 element, chosen orthogonal to $h$.
So $B=\{ 1,h,c,h^2,h^3,h^4 \}$
with $c\cdot h=0$, $\int h^4\ne 0$, $\int c^2\ne 0$.
We have $K=-4h$, $K=-3h$ in the cases of the two respective varieties;
set $\kappa=4$, $\kappa=3$ accordingly.

The relations in curve class $\beta$ never involve the
number $N(\beta; \kappa\beta+1,0,0,0)$.
We recover this exceptional number from a particular degree $\beta+1$
relation in which it appears in a quadratic term, for which, to be
able to solve, we must add the hypothesis that
$N(1;0,0,1,1)\ne0$ (resp.\ $N(1;1,0,0,1)\ne0$) when $\kappa=4$
(resp.\ $\kappa=3$).

For each $\beta$, there are 10 numbers $N(\beta;d)$ which are not
of the form $N(\beta;t,u,v,w)$ with $u+v+w\ge3$; these are the
numbers unreachable by the proof of strong reconstruction
applied to the subring of $A$ generated by $h$.
We must show how to solve for 9 of these (all except
$N(\beta; \kappa\beta+1,0,0,0)$) as well as the
leftover degree $\beta-1$ number.
Table \ref{gtwofour} outlines how to do this.

\begin{table}
\begin{center}
\begin{tabular}{lcclr}
\smallskip (a) & $N(\beta; \kappa\beta-5,0,0,2)$ & by &
$\dfeyn {h^4}h{h^3}c{(\beta;\kappa\beta-6,0,0,0)}$ & $(\beta\ge2)$ \\
\smallskip (b) & $\fra{\displaystyle N(\beta-1;\hfill}
{\displaystyle\ \kappa\beta-\kappa+1,0,0,0)}$ & by &
$\dfeyn ch{h^4}c{(\beta;\kappa\beta-5,0,0,0)}$ & $(\beta\ge2)$ \\
\smallskip (c) & $\fra{\displaystyle N(\beta; t,u,v,w)\hfill}
{\displaystyle{\rm with\ }u+v+w\ge3}$ & by &
$\dfeyn {\langle h \rangle}{\langle h \rangle}{\langle h \rangle}
{\langle h \rangle}{(\beta;d)}$ & \\
\smallskip (d) & $N(\beta; \kappa\beta-4,0,1,1)$ & by &
$\dfeyn {h^4}h{h^2}c{(\beta;\kappa\beta-5,0,0,0)}$ & $(\beta\ge2)$ \\
\smallskip (e) & $N(\beta; \kappa\beta-3,1,0,1)$ & by &
$\dfeyn {h^4}hhc{(\beta;\kappa\beta-4,0,0,0)}$ & $(\kappa+\beta\ge5)$ \\
\smallskip (f) & $N(\beta; \kappa\beta-3,0,2,0)$ & by &
$\dfeyn {h^3}h{h^2}c{(\beta;\kappa\beta-4,0,0,0)}$ & $(\kappa+\beta\ge5)$ \\
\smallskip (g) & $N(\beta; \kappa\beta-2,0,0,1)$ & by &
$\dfeyn {h^3}hcc{(\beta;\kappa\beta-4,0,0,0)}$ & $(\kappa+\beta\ge5)$ \\
\smallskip (h) & $N(\beta; \kappa\beta-2,1,1,0)$ & by &
$\dfeyn {h^3}hhc{(\beta;\kappa\beta-3,0,0,0)}$ & \\
\smallskip (i) & $N(\beta; \kappa\beta-1,0,1,0)$ & by &
$\dfeyn {h^2}hcc{(\beta;\kappa\beta-3,0,0,0)}$ & \\
\smallskip (j) & $N(\beta; \kappa\beta-1,2,0,0)$ & by &
$\dfeyn {h^2}hhc{(\beta;\kappa\beta-2,0,0,0)}$ & \\
\smallskip (k) & $N(\beta; \kappa\beta,1,0,0)$ & by &
$\dfeyn hhcc{(\beta;\kappa\beta-2,0,0,0)}$ &
\end{tabular}
\end{center}
\caption{Order, within the outer induction, for the proof
of strong reconstruction for $G(2,4)$/$\Sym^2\bP^2$.}
\label{gtwofour}
\end{table}

We induct first on curve class $\beta$.
The induction hypothesis consists of
all relations and
all numbers in degrees less than $\beta-1$, and all
relations and all numbers except
$N(\beta-1; \kappa\beta-\kappa+1,0,0,0)$ in degree $\beta-1$.
The relations indicated in entries (a), (b) of Table \ref{gtwofour}
give us two new numbers, including $N(\beta-1; \kappa\beta-\kappa+1,0,0,0)$.
Thus from now on we assume all relations and
all numbers in degrees less than or equal to $\beta-1$.

Now, inductively on $d$ via the partial ordering
$d'=(t',u',v',w')\prec d=(t,u,v,w)\Leftrightarrow
|d'| < |d|$ or $|d'|=|d|$, $t'>t$ or
$|d'|=|d|$, $t'=t$, $u'+2v'+3w' > u+2v+3w$, we establish
the relations indicated in (c)--(k) of the table, as applicable
to the current degree.
Step (c) is the inner induction of the proof
of the strong reconstruction theorem,
applied to the subring of $A$ generated by $h$.

Finally, it follows from \fsr{} (this takes a bit of checking)
that in any degree, the (31 out of 55 total) relations indicated
in Table \ref{gtwofour} imply all the relations.
The 21 encoded by step (c) follow by the proof of strong reconstruction,
so we are reduced to establishing the remaining 10.
When $u=v=w=0$, we are done by Table \ref{gtwofour} (each remaining relation
solves for an unknown number).
Otherwise, since each of the Feynman diagrams indicated
in (a), (b), (d)--(k) of the table has a diagonal
containing $c$ and $h$,
an application of \tsr{} reduces us to relations obtained
in previous degrees (via the partial ordering above).

Summarizing, in a neighborhood of the geometric solution,
the solution space to WDVV consists of only the expected rescalings
for $X=G(2,4)$, but has a dependence on two extra free parameters
in case $X=\Sym^2\bP^2$
(where the geometric solution comes from viewing $X$ as a
global quotient of a homogeneous variety by a finite group).
When $N(1;0,0,1,1)=0$ (resp.\ $N(1;1,0,0,1)=0$) the situation is
considerably more complicated.
For instance, for $X=G(2,4)$, setting $N(1;5,0,0,0)=1$
and all other $N(1;t,u,v,w)=0$ yields a consistent solution through $\beta=4$
(with $N(\beta;4\beta+1,0,0,0)$ indeterminate for $\beta=2,3,4$),
but fails to satisfy the relations consistently starting in degree 5.

{\bf Example 4.} $X=G(2,5)$.
The same technique as that of Example 3 establishes strong
reconstruction for $G(2,5)$.
The statement of strong reconstruction dictates 3 starting numbers
and no relationships, but as we shall see, we need to assume
one relationship.
As in the previous example, we will also have to make a genericity assumption.
We shall find that the geometric solution to WDVV lives in
a two-parameter family of solutions.

\begin{table}
\begin{center}
\begin{tabular}{c|cccccccc}
 & $t_1$ & $t_2$ & $t_3$ & $t_4$ & $t_5$ & $t_6$ & $t_7$ & $t_8$ \\ \hline
$t_1$ & $t_2$ & $t_4$ & $t_5$ & $t_6$ & $(1/3)t_7$ & $5t_8$ & $0$ & $t_9$ \\
$t_2$ & $t_4$ & $t_6$ & $(1/3)t_7$ & $5t_8$ & $0$ & $5t_9$ & $0$ & $0$ \\
$t_3$ & $t_5$ & $(1/3)t_7$ & $t_6-(11/3)t_7$ & $0$ & $5t_8$ &
$0$ & $15t_9$ & $0$ \\
$t_4$ & $t_6$ & $5t_8$ & $0$ & $5t_9$ & $0$ & $0$ & $0$ & $0$ \\
$t_5$ & $(1/3)t_7$ & $0$ & $5t_8$ & $0$ & $5t_9$ & $0$ & $0$ & $0$ \\
$t_6$ & $5t_8$ & $5t_9$ & $0$ & $0$ & $0$ & $0$ & $0$ & $0$ \\
$t_7$ & $0$ & $0$ & $15t_9$ & $0$ & $0$ & $0$ & $0$ & $0$ \\
$t_8$ & $t_9$ & $0$ & $0$ & $0$ & $0$ & $0$ & $0$ & $0$
\end{tabular}
\end{center}
\caption{Multiplication table for $A^*_\Q G(2,5)$ with respect to the
basis of T. Graber.}
\label{multtable}
\end{table}

Thinking of $X$ as the space of rank 2 quotients of
$\C^5$, let $Q$ be the universal quotient bundle and $c_i=c_i(Q)$.
We use the following basis for $A^*_\Q X$, suggested by T. Graber,
\begin{center}
\begin{tabular}{llll}
codim 0:&& $t_0=1$ \\
codim 1:&& $t_1=c_1$ \\
codim 2:&& $t_2=c_1^2$ & $t_3=2c_1^2-5c_2$ \\
codim 3:&& $t_4=c_1^3$ & $t_5=2c_1^3-5c_1c_2$ \\
codim 4:&& $t_6=c_1^4$ & $t_7=c_1^4-5c_2^2$ \\
codim 5:&& $t_8=c_1c_2^2$ \\
codim 6:&& \multicolumn{2}{l}{ $t_9=c_2^3$ (point class)} \\
\end{tabular}
\end{center}
with multiplication table given in Table \ref{multtable}.

We denote a typical unknown by
$N(\beta;d_2,d_3,d_4,d_6,d_8,d_9,d_3,d_5,d_7)$ (note special order).
Inductively on degree $\beta$, we show that degree $\beta$
relations solve consistently for all but 8 degree $\beta$ numbers, plus
the 5 linear expressions shown in Table \ref{linexpr}.
The 8 exceptions are the 7 numbers appearing in Table \ref{linexpr}
as well as $N(\beta;0,0,0,0,0,5\beta+3,0,0)$.

The genericity assumption is $N(1;0,0,0,0,1,0,0,1) \ne 0$.
Given the induction hypothesis, we solve for the remaining degree
$(\beta-1)$ numbers according to Table \ref{gtwofivea},
first by the path shown
with $(0,0,0,0,0,-4,2,0)$ added to all degrees, then by the
path shown with $(0,0,0,0,0,-2,1,0)$ added to all degrees,
and then by the path as shown.
Only for $\beta=2$, during the first pass, we must substitute
$\dfeyn {t_9}{t_1}{t_6}{t_5}{(\beta;0,0,0,0,0,5\beta-10,1,0)}$
for step (c).

Once we have all the numbers in degree $\beta-1$, we then
induct on $d$ with respect to the partial ordering
$d'\prec d\Leftrightarrow$
\begin{itemize}
\item[(i)] $|d'|<|d|$, or
\item[(ii)] $|d'|=|d|$ and
$d'_2+d'_4+d'_6+d'_8+d'_9<d_2+d_4+d_6+d_8+d_9$, or
\item[(iii)] $|d'|=|d|$ and
$d'_2+d'_4+d'_6+d'_8+d'_9=d_2+d_4+d_6+d_8+d_9$, but
$d'_2+2d'_4+3d'_6+4d'_8+5d'_9>d_2+2d_4+3d_6+4d_8+5d_9$, or
\item[(iv)] $|d'|=|d|$, $(d'_2,d'_4,d'_6,d'_8,d'_9)=(d_2,d_4,d_6,d_8,d_9)$,
and $d'_3+2d'_5+3d'_7<d_3+2d_5+3d_7$.
\end{itemize}
For each $d$, we perform the inner induction of the proof of
strong reconstruction to obtain all
relations involving only powers of $t_1$.
Next, we obtain all of
$$
\begin{tabular}{lll}
\smallskip $\fra{\displaystyle N(\beta;0,0,0,0,2,u,v,w)}
{\displaystyle \hfill u+2v+3w=5\beta-7}$ & by &
$\dfeyn {t_9}{t_1}{t_8}{t_7}{(\beta;0,0,0,0,0,u,v,w-1)}$ \\
\smallskip & or & $\dfeyn {t_9}{t_1}{t_8}{t_5}{(\beta;0,0,0,0,0,u,v-1,w)}$ \\
\smallskip & or & $\dfeyn {t_9}{t_1}{t_8}{t_3}{(\beta;0,0,0,0,0,u-1,v,w)}$ \\
\smallskip $\fra{\displaystyle N(\beta;0,0,0,1,1,u,v,w)}
{\displaystyle \hfill u+2v+3w=5\beta-6}$ & by &
$\dfeyn {t_9}{t_1}{t_6}{t_7}{(\beta;0,0,0,0,0,u,v,w-1)}$ \\
& or & etc. \\
$\ldots$ \\
$\fra{\displaystyle N(\beta;2,0,0,0,0,u,v,w)}
{\displaystyle \hfill u+2v+3w=5\beta+1}$ & by &
$\dfeyn {t_2}{t_1}{t_1}{t_7}{(\beta;0,0,0,0,0,u,v,w-1)}$ etc.\\
\end{tabular}
$$
coming from relations in degree $(\beta;d)$.
The exception to be noted occurs in attempting to solve for
$N(\beta;2,0,0,0,0,5\beta+1,0,0)$: we get a value
for $(L5)$ of Table \ref{linexpr} rather than a single $N$.

\begin{table}
$$\begin{array}{ll}
\medskip (L1) & N(\beta;0,0,0,0,0,5\beta,0,1) -
3 N(\beta;0,0,0,0,0,5\beta-1,2,0) \\
(L2) & N(\beta;1,0,0,0,0,5\beta,1,0) - N(\beta;0,1,0,0,0,5\beta+1,0,0) \\
\medskip & \hspace{130pt} {} + \beta N(\beta;0,0,0,0,0,5\beta-1,2,0) \\
\medskip (L3) & N(\beta;1,0,0,0,0,5\beta,1,0) -
2 \beta N(\beta;0,0,0,0,0,5\beta-1,2,0) \\
(L4) & N(\beta;1,0,0,0,0,5\beta+2,0,0) -
2\beta N(\beta;0,0,0,0,0,5\beta+1,1,0) \\
\medskip & \hspace{130pt} {} -
11 \beta^2 N(\beta;0,0,0,0,0,5\beta-1,2,0) \\
(L5) & N(\beta;2,0,0,0,0,5\beta+1,0,0) -
4 \beta^2 N(\beta;0,0,0,0,0,5\beta-1,2,0)
\end{array}$$
\caption{The linear expressions obtained by degree $\beta$ relations.}
\label{linexpr}
\end{table}

Still in a particular degree, we obtain
$$
\begin{tabular}{cll}
\smallskip $\fra{\displaystyle N(\beta;0,0,0,0,1,u,v,w)}
{\displaystyle \hfill u+2v+3w=5\beta-2}$ & by &
$\dfeyn{t_8}{t_1}{t_7}{t_7}{(\beta;0,0,0,0,0,u,v,w-2)}$ \\
 & or & $\dfeyn{t_8}{t_1}{t_7}{t_5}{(\beta;0,0,0,0,0,u,v-1,w-1)}$ \\
 & etc.
\end{tabular}
$$
with exceptions noted below:
$$
\begin{tabular}{ccl}
\smallskip \hspace{80pt}$(L2)$
& by & $\dfeyn{t_2}{t_1}{t_3}{t_3}{(\beta;0,0,0,0,0,5\beta-1,0,0)}$ \\
\smallskip \hspace{80pt}$(L3)$
& by & $\dfeyn{t_1}{t_1}{t_5}{t_3}{(\beta;0,0,0,0,0,5\beta-1,0,0)}$ \\
\hspace{80pt}$(L4)$
& by & $\dfeyn{t_1}{t_1}{t_3}{t_3}{(\beta;0,0,0,0,0,5\beta,0,0)}$
\end{tabular}
$$

\begin{table}
\begin{center}
\begin{tabular}{lccl}
\smallskip (a) & $N(\beta; 0,0,1,0,1,5\beta-5,0,0)$ & by &
$\dfeyn {t_9}{t_1}{t_4}{t_3}{(\beta;0,0,0,0,0,5\beta-6,0,0)}$ \\
\smallskip (b) & $N(\beta; 0,0,0,2,0,5\beta-5,0,0)$ & by &
$\dfeyn {t_6}{t_1}{t_8}{t_3}{(\beta;0,0,0,0,0,5\beta-6,0,0)}$ \\
\smallskip (c) & $N(\beta; 0,0,0,1,1,5\beta-6,0,0)$ & by &
$\dfeyn {t_9}{t_1}{t_6}{t_3}{(\beta;0,0,0,0,0,5\beta-7,0,0)}$ \\
\smallskip (d) & $N(\beta; 0,0,0,0,1,5\beta-5,0,1)$ & by &
$\dfeyn {t_8}{t_1}{t_7}{t_3}{(\beta;0,0,0,0,0,5\beta-6,0,0)}$ \\
\smallskip (e) & $N(\beta; 0,0,0,1,0,5\beta-4,0,1)$ & by &
$\dfeyn {t_6}{t_1}{t_7}{t_3}{(\beta;0,0,0,0,0,5\beta-5,0,0)}$ \\
\smallskip (f) & $N(\beta; 0,0,0,0,1,5\beta-4,1,0)$ & by &
$\dfeyn {t_8}{t_1}{t_5}{t_3}{(\beta;0,0,0,0,0,5\beta-5,0,0)}$ \\
\smallskip (g) & $N(\beta{-}1; 0,0,0,0,0,5\beta-2,0,0)$ & by &
$\dfeyn {t_3}{t_1}{t_9}{t_3}{(\beta;0,0,0,0,0,5\beta-5,0,0)}$
\end{tabular}
\end{center}
\caption{Path to the remaining degree $(\beta-1)$ numbers for $G(2,5)$.}
\label{gtwofivea}
\end{table}

Lastly, we have numbers of the form $N(\beta;0,0,0,0,0,u,v,w)$
and the relations that produce these:
$$
\begin{tabular}{cccclcccc}
relation & \multicolumn{3}{c}{for cases} & & relation &
\multicolumn{3}{c}{for cases} \\
\smallskip $\feyn{t_5}{t_1}{t_7}{t_7}$ &    & & $w\ge3$  & &
$\feyn{t_3}{t_1}{t_7}{t_3}$ &  $u\ge1$ & $v\ge1$ & $w\ge1$ \\
\smallskip $\feyn{t_5}{t_1}{t_7}{t_5}$ &   & $v\ge1$ & $w\ge2$ & &
$\feyn{t_5}{t_1}{t_3}{t_5}$ &   &  $v\ge3$ \\
\smallskip $\feyn{t_5}{t_1}{t_7}{t_3}$ &  $u\ge1$ & & $w\ge2$ & &
$\feyn{t_3}{t_1}{t_5}{t_3}$ & \multicolumn{3}{r}{(only to get $(L1)$)} \\
$\feyn{t_3}{t_1}{t_7}{t_5}$ &  & $v\ge2$ & $w\ge1$
\end{tabular}
$$

Finally, we obtain $\feyn {t_3}{t_1}{t_9}{t_3}$,
$\feyn {t_5}{t_1}{t_9}{t_3}$, and
$\feyn {t_7}{t_1}{t_9}{t_3}$.
When $(d_2,d_4,d_6,d_8,d_9)=(0,0,0,0,0)$ and $d_5=d_7=0$,
these are established by the three passes
through the path in Table \ref{gtwofivea}.
For the details see Table \ref{gtwofiveb}.
Note that each pass through Table \ref{gtwofivea} is used to solve for
an unknown when $\beta\ge2$; when $\beta=1$ we
actually get a condition on starting data, described below.

\begin{table}
\begin{center}
\begin{tabular}{crcl}
Assumptions & \multicolumn{3}{c}{Implications} \\
$\dfeyn{t_8}{t_1}{t_7}{t_3}{(5\beta-8,1,0)}$ &
$\dfeyn{t_3}{t_1}{t_9}{t_3}{(5\beta-9,2,0)} + \mbox{\tsr}$
& $\Rightarrow$ & $\dfeyn{t_3}{t_1}{t_9}{t_5}{(5\beta-8,1,0)}$ \\
$\dfeyn{t_5}{t_1}{t_8}{t_5}{(5\beta-8,1,0)}$ &
$\fiveblu{t_1}{t_6}{t_1}{t_7}{t_3}$
& $\Rightarrow$ & $\dfeyn{t_8}{t_1}{t_3}{t_7}{(5\beta-8,1,0)}$ \\
$\dfeyn{t_6}{t_1}{t_7}{t_5}{(5\beta-8,1,0)}$ &
$\fivecru{t_5}{t_1}{t_8}{t_1}{t_3}$
& $\Rightarrow$ & $\dfeyn{t_5}{t_1}{t_9}{t_3}{(5\beta-8,1,0)}$ \\
$\dfeyn{t_8}{t_1}{t_7}{t_5}{(5\beta-7,0,0)}$ &
\tsr & $\Rightarrow$ & $\dfeyn{t_5}{t_1}{t_9}{t_5}{(5\beta-7,0,0)}$ \\
& $\fivecru{t_1}{t_7}{t_1}{t_8}{t_3}$
& $\Rightarrow$ & $\dfeyn{t_1}{t_7}{t_9}{t_3}{(5\beta-7,0,0)}$ \\
& $\fiveblu{t_1}{t_5}{t_1}{t_9}{t_3}$
& $\Rightarrow$ & $\dfeyn{t_7}{t_1}{t_9}{t_3}{(5\beta-7,0,0)}$
\end{tabular}
\end{center}
\caption{How $\dfeyn {t_3}{t_1}{t_9}{t_3}{(5\beta-9,2,0)}\Rightarrow
\dfeyn {t_7}{t_1}{t_9}{t_3}{(5\beta-7,0,0)}$
follows from \tsr{} and \fsr{}.
We use $(u,v,w)$ as shorthand for degree
$(\beta;0,0,0,0,0,u,v,w)$.
Each relation listed as an assumption
solves for an unknown number and for that reason is satisfied.
Note that $\dfeyn {t_3}{t_1}{t_9}{t_3}{(5\beta-7,1,0)}\Rightarrow
\dfeyn {t_5}{t_1}{t_9}{t_3}{(5\beta-6,0,0)}$
is just the first three steps above.}
\label{gtwofiveb}
\end{table}

When $(d_2,d_4,d_6,d_8,d_9)=(0,0,0,0,0)$ but $d_5\ne0$ or $d_7\ne0$, then
because of the induction order, some of the (85 total) relations listed
as determining numbers will determine numbers that have already
been solved for.
But in each such case, \tsr{} allows us to deduce the relation in question.
Finally, when $(d_2,d_4,d_6,d_8,d_9)\ne(0,0,0,0,0)$
then all 85 relations follow by \tsr{} just as in Example 3.

As in Example 3, we now note that the 85 relations in the above lists
plus the 120 only involving powers of $t_1$ imply the remaining
461 relations by \fsr.
Because of the number of relations involved, the author used a computer
to complete this verification.

Relation 
$\dfeyn {t_3}{t_1}{t_9}{t_3}{(\beta;0,0,0,0,0,5\beta-5,0,0)}$
from (i) in Table \ref{gtwofivea} solves for an unknown number
whenever $\beta\ge2$,
but still needs to be verified when $\beta=1$.
This is what imposes the one condition on starting data.
When written out for $\beta=1$, (e)--(g) of Table \ref{gtwofivea} translate as
\begin{eqnarray*}
\lefteqn{11 N(1;0,0,0,0,1,0,0,1)=} \hspace{30pt} \\
& & 6 N(1;0,0,1,0,1,0,0,0)+15 N(1;0,0,0,2,0,0,0,0).
\end{eqnarray*}

\section{Appendix: geometric groundwork}
\label{geoground}

Let $X$ be a complex projective manifold, and for simplicity,
assume $X$ has homology only in even dimensions.
To avoid writing doubled indices, set $A_dX=H_{2d}(X,\Z)$ and
$A^dX=H^{2d}(X,\Z)$.
Set $A^d_\Q X=A^dX\otimes\Q$.
For $\beta\in A_1X$, we denote by $M_{0,n}(X,\beta)$
the moduli space of $n$-pointed rational curves $f\colon \bP^1\to X$
such that $f_*[\bP^1]=\beta$.
The Kontsevich space
$\overline{M}_{0,n}(X,\beta)$ is a compactification of $M_{0,n}(X,\beta)$
whose points correspond to $n$-pointed trees of $\bP^1$'s mapping
into $X$ satisfying a stability hypothesis.
Evaluation maps $\rho_i\colon\overline{M}_{0,n}(X,\beta)\to X$ send
a given map $f$ to the image under $f$ of the $i^{\rm th}$ marked point.

If $X$ is a homogeneous variety (a quotient of a complex reductive Lie
group by a parabolic subgroup), then the
{\em tree-level system of Gromov-Witten numbers} is the system of maps
$I_\beta\colon \bigoplus \Sym^n A^*_\Q X\to\Q$
for each $\beta\in A_1X$
given by the formula
\begin{equation}
\label{gwformula}
I_\beta(\gamma_1\cdots\gamma_n)=\int_{\overline{M}_{0,n}(X,\beta)}
\rho_1^*(\gamma_1)\smallcup\cdots\smallcup\rho_n^*(\gamma_n).
\end{equation}
With the notation of Problem \ref{mainprob},
the {\em geometric solution} to WDVV is given by
$N(\beta;d)=I_\beta(\bigotimes (T_{\tau_j})^{\otimes d_j})$.
These numbers have enumerative significance: $N(\beta;d)$ is the number of
rational curves on $X$ in
homology class $\beta$ meeting $d_j$ general translates
of a cycle Poincar\'e dual to $T_{\tau_j}$, for each $j$.

For general $X$, there is still a geometric solution to WDVV,
given by formula (\ref{gwformula}) but with the
integration performed over a virtual fundamental cycle,
although the enumerative significance of these numbers is
less clear.
It is still the case for any $\beta\ne0$ that
$I_\beta=0$ unless $\int_\beta\omega>0$
for every ample divisor $\omega$.

There exist forgetful maps forgetting $X$ and forgetting any subset
of the marked points, so in particular to any
$\{i,j,k,l\}\subset\{1,\ldots,n\}$ there corresponds a map
$\overline{M}_{0,n}(X,\beta)\to\overline{M}_{0,\{i,j,k,l\}}$.
Pulling back rational equivalences on $\overline{M}_{0,4}\cong\bP^1$
leads to the associativity relations
that the Gromov-Witten numbers must satisfy.

Since the differential equation (\ref{assrel})
is invariant under translations
$y_{\sigma_i}\mapsto y_{\sigma_i}+\alpha_i$, there is
an $r$-dimensional family of rescalings acting on the
geometric solution to WDVV, where $r=\rk A^1X$.
Much of the contents of this paper was motivated by
a search for solutions to WDVV besides the geometric
solution and its translates under these rescalings.

The formal computation of section \ref{fivesi} is motivated
by considering a typical
component of the boundary of one of the Kontsevich spaces
$\overline{M}_{0,n}(X,\beta)$.
Say $A\cup B=\{1,\ldots,n\}$, $A\cap B=\emptyset$, $|A|\ge2$, $|B|\ge2$,
and $\beta_1+\beta_2=\beta$.
Then there is a component $D(A,B;\beta_1,\beta_2)$ of the boundary
of $\overline{M}_{0,n}(X,\beta)$, which
fits into a fiber diagram (see \cite{c})
\begin{equation*}
\begin{CD}
 D(A,B;\beta_1,\beta_2) @>>>  X^n\times X \\
@VVV @V{\delta}VV \\
\overline{M}_{0,A\cup\{{*}\}}(X,\beta_1)\times
\overline{M}_{0,B\cup\{{*}\}}(X,\beta_2) @>{\rho}>> X^n\times X\times X \\
\end{CD}
\end{equation*}
where $\delta$ is given by the diagonal embedding of $X$ in $X\times X$.

Given cohomology classes $\gamma_1,\ldots,\gamma_n,\xi\in A^*X$,
we have
\begin{eqnarray*}
\delta_*(\gamma_1\times\cdots\times\gamma_n\times\xi) &=&
\sum_{e,f} g^{ef} \gamma_1\times\cdots\times\gamma_n
\times(\xi\smallcup T_e)\times T_f \\
&=& \sum_{e,f} g^{ef} \gamma_1\times\cdots\times\gamma_n
\times T_e\times (\xi\smallcup T_f)
\end{eqnarray*}
corresponding to two ways of writing the class $\xi$ on
the diagonal.
Pulling back by $\rho$ and integrating gives the identity
\begin{eqnarray*}
\lefteqn{ \sum_{e,f} g^{ef} I_{\beta_1}\bigl(\prod_{a\in A}
\gamma_a \cdot (\xi\smallcup T_e)\bigr)
I_{\beta_2}\bigl(\prod_{b\in B} \gamma_b\cdot T_f\bigr) = }
\hspace{20pt} \\
& & \sum_{e,f} g^{ef} I_{\beta_1}
\bigl(\prod_{a\in A}\gamma_a \cdot T_e\bigr)
I_{\beta_2}\bigl(\prod_{b\in B}\gamma_b\cdot
(\xi\smallcup T_f)\bigr).
\end{eqnarray*}
If $\gamma_1=T_i$, $\gamma_2=T_j$, $\gamma_3=T_k$, $\gamma_4=T_l$,
and $\xi=T_m$, and if we sum over partitions $(\beta_1,\beta_2,A,B)$
such that $A$ contains 1 and 2 and $B$ contains 3 and 4,
we get a special case of (\ref{mdiag}), namely the case where
the potential function $\ga$ is given by the geometric solution
to WDVV.

\noindent
Department of Mathematics\\
University of Chicago\\
5734 S. University Avenue\\
Chicago, IL 60637\\
kresch\atsign{}math.uchicago.edu
\end{document}